\documentclass{amsart}
\usepackage{palatino,vmargin}
\newtheorem{thm}{Theorem}[section]

\newtheorem{lem}[thm]{Lemma}
\newtheorem{prop}[thm]{Proposition}
\theoremstyle{definition}

\theoremstyle{remark}
\newtheorem{rem}[thm]{Remark}

\numberwithin{equation}{section}
\newcommand{\mylabel}[1]{{\label{#1}}}%{{\fbox{#1}}}}
\newcommand{\inta}[1]{\mathop{{\displaystyle\int\kern-1.2em\raise2.5pt\hbox to 10pt{\leaders\hrule width 8pt height1pt\hfill}}}_{#1}}
\begin{document}
\title{ The mixed problem for harmonic functions in polyhedra of $\mathbb{R}^3$}%
\author{Moises Venouziou and Gregory C. Verchota}%
\thanks{The second author gratefully acknowledges partial support provided by
the National Science Foundation  through award DMS-0401159\\
\textit{E-mail address:}
gverchot@syr.edu\;\;\;\;\;\;\;\;\;\;\;\;\;\;\;\;\;\;\;\;\;\;\;\;\;\;\;\;\;\;}
\address{215 Carnegie \\ Syracuse University\\ Syracuse NY 13244}%
 \email{gverchot@syr.edu}
\date{\today}%
\begin{abstract}
R. M. Brown's theorem on mixed Dirichlet and Neumann boundary
conditions is extended in two ways for the special case of
polyhedral domains.  A (1) more general partition of the boundary
into Dirichlet and Neumann sets is used on (2) manifold boundaries
that are not locally given as the graphs of functions.  Examples are
constructed to illustrate  necessity and other implications of the
geometric hypotheses.
\end{abstract}
\maketitle
\section{Introduction}\mylabel{Introduction}

In \cite{Bro94} R. M. Brown initiated a study of the \textit{mixed
boundary value problem} for harmonic functions in \textit{creased
Lipschitz domains} $\Omega$ with data in the Lebesgue and Sobolev
spaces $L^2(\partial\Omega)$ and $W^{1,2}(\partial\Omega)$ (with
respect to surface measure $ds$) taken in the strong pointwise sense
of nontangential convergence.

At the end of his article Brown poses a question concerning a
certain topologic-geometric difficulty not included in his solution:
Can the mixed problem be solved in the (infinite) pyramid of
$\mathbb{R}^3$,  $|X_1|+|X_2|< X_3$,  when Neumann and Dirichlet
data are chosen to alternate on the faces?  In this article we avoid
the geometric difficulties of what can be called Lipschitz faces or
facets and provide answers in the case of compact polyhedral domains
of $\mathbb{R}^3$.  Some other recent approaches to the mixed
problem for second order operators and systems in polyhedra can be
found in \cite{MR07} \cite{MR06} \cite{MR05} \cite{MR04} \cite{MR03}
\cite{MR02} and \cite{Dau92}.

Consider a compact polyhedron of $\mathbb{R}^3$ with the property
that its interior $\Omega$ is connected.  $\Omega$ will be termed a
\textit{compact polyhedral domain}.  Suppose its boundary
$\partial\Omega$ is a connected $2$-manifold.  Such a domain
$\Omega$ need not be a Lipschitz domain.  Partition the boundary of
$\Omega$ into two disjoint sets $N$ and $D$, for Neumann data and
Dirichlet data respectively, so that the following is satisfied.

\begin{multline}\mylabel{partition}
\mbox{(i) $N$ is the union of a number (possibly \textit{zero}) of
closed faces of $\partial\Omega$.}\\  \mbox{(ii)
$D=\partial\Omega\setminus N$ is nonempty.\hskip3.42in}\\
\mbox{(iii) Whenever a face of $N$ and a face of $D$ share a
$1$-dimensional edge as boundary, the}\\\mbox{ dihedral angle
measured in $\Omega$ between the two faces is \textit{less} than
$\pi$.\hskip1.2in}
\end{multline}

The $L^2$-\textit{polyhedral mixed problem} for harmonic functions
is
\begin{multline}\mylabel{mixedproblem}
\mbox{  Given $f\in W^{1,2}(\partial\Omega)$ and $g\in L^2(N)$ show
there exists a  solution to $\triangle u=0$ in $\Omega$ such
that}\\\mbox{ (i) $u\rightarrow^{n.t.}f$  $a.e.$ on $D$.\hskip4in}
\\\mbox{ (ii) $\partial_\nu u\rightarrow^{n.t.}g$  $a.e.$ on
$N$.\hskip3.9in}\\\mbox{ (iii) $\nabla u^*\in
L^2(\partial\Omega)$.\hskip4.4in}
\end{multline}
Here $\nabla u^*$ is the \textit{nontangential maximal function} of
the \textit{gradient} of $u$.  Generally for a function $w$ defined
in a domain $G$
$$w^*(P)=\sup_{X\in\Gamma(P)}|w(X)|, P\in \partial G.$$
For a choice of $\alpha>0$ \textit{nontangential approach regions}
for each $P\in \partial G$ are defined by

\begin{equation}\mylabel{approach}\Gamma(P)=\{X\in G: |X-P|<(1+\alpha) dist(X,\partial G)\}\end{equation}
Varying the choice of $\alpha$ yields nontangential maximal
functions with comparable $L^p(\partial G)$ norms $1<p\leq\infty$ by
an application of the Hardy-Littlewood maximal function.  Therefore
$\alpha$ is suppressed.  In  general when $w^*\in L^p(\partial G)$
is written it is understood that the nontangential maximal function
is with respect to cones determined by the domain $G$.  The
\textit{outer unit normal} vector to $\Omega$ (or a domain $G$) is
denoted $\nu=\nu_P$ for $a.e. P\in \partial \Omega$ and the limit of
(ii) is understood as $$\lim_{ \Gamma(P)\ni X\rightarrow P}
\nu_P\cdot \nabla u (X)=g(P)$$ and similarly for (i).

A consequence of solving \eqref{mixedproblem} is that the gradient
of the solution has well defined nontangential limits at the
boundary $a.e.$

In addition, as Brown points out, solving the mixed problem yields
\textit{extension operators} $W^{1,2}(D)\rightarrow W^{1,2}(\partial
\Omega)$ by $f\mapsto u|_{\partial \Omega}$ where $u$ is a solution
to the mixed problem with $u|_D=f$.  Consequently \textit{problem
\eqref{mixedproblem} cannot be solved for all $f\in W^{1,2}(D)$ when
$D$ and $N$ are defined as on the boundary of the pyramid}.  For
example, since the pyramid is Lipschitz at the origin so that
Sobolev functions on its boundary project to Sobolev functions on
the plane, solving \eqref{mixedproblem} implies that a local
$W^{1,2}$ function exists in $\mathbb{R}^2$ that is identically $1$
in the first quadrant and identically $zero$ in the third.  Such a
function necessarily restricts to a local $W^{\frac{1}{2},2}$
function on any straight line  through the origin. But a step
function is not locally in $W^{\frac{1}{2},2}(\mathbb{R})$.  The
boundary domain $D$ (and its projection) do not satisfy the
\textit{segment property} commonly invoked to show the two Sobolev
spaces $H_1(D)$ and $W^{1,2}(D)$ equal \cite{Agm65} \cite{GT83}.

The admissible Sobolev functions on $D$ must then be those that have
extensions to $W^{1,2}(\partial \Omega)$.  Or equivalently, the
admissible Sobolev functions on $D$ are the restrictions of
$W^{1,2}(\partial \Omega)$ functions. We introduce the following
\textit{norm} on the \textit{space of restrictions of
$W^{1,2}(\partial \Omega)$ functions $f$ to $D$}
$$\|f\|^2_{D}=\inf_{\widetilde{f}|_D=f}\int _{\partial \Omega}\widetilde{f}^2+
|\nabla_t \widetilde{f}|^2 ds$$ Here $\widetilde{f}$ denotes all
$W^{1,2}(\partial \Omega)$ functions that restrict to $f$ on $D$,
and $\nabla_t$ denotes the tangential gradient.  That this is a norm
follows by arguments such as:  Given $f\in W^{1,2}(\partial \Omega)$
and a real number $a$,  the functions $a\widetilde{f}$ form a subset
of all extensions $\widetilde{af}$ of $(af)|_D$ so that
$\|af\|_D\leq |a|\|f\|_D$, and thus likewise $\|f\|_D\leq
|a|^{-1}\|af\|_D$ when $a \neq 0$.

This normed space is complete by using the standard completeness
proof for Lebesgue spaces:  Given a Cauchy sequence $\{f_j\}$ let
$j_k$ be such that $\|f_i-f_j\|_D<2^{-k}$ for all $i,j\geq j_k$ and
define $g_k=f_{j_{k+1}}-f_{j_{k}}$.  Then there exists an extension
$\widetilde{g_k}$ such that $\|\widetilde{g_k}\|_{W^{1,2}(\partial
\Omega)}< 2^{-k}$.  Extensions of $f_{j_{n+1}}$ may then be defined
by $\widetilde{f_{j_{1}}}+\sum_{k=1}^n \widetilde{g_k}$ Cauchy in
$W^{1,2}(\partial \Omega)$.  Completeness will follow.  The Banach
space of restrictions to $D$ is undoubtedly the generally smaller
Sobolev space $H_1(D)$ (e.g. \cite{Fol95} p. 220), but this will not
be pursued further.

A homogeneous Sobolev semi-norm on $D$ is defined by

\begin{equation}\mylabel{seminorm}
\|f\|^2_{D^o}=\inf_{\widetilde{f}|_D=f}\int _{\partial \Omega}
|\nabla_t \widetilde{f}|^2 ds
\end{equation}

When $\partial \Omega$ is connected the following \textit{scale
invariant} theorem is established in the Section \ref{proof}.

\begin{thm}\mylabel{mixed}
Let $\Omega\subset \mathbb{R}^3$ be a compact polyhedral domain with
connected $2$-manifold boundary $\partial\Omega=D\cup N$ satisfying
the conditions \eqref{partition}.  Then given $f\in
W^{1,2}(\partial\Omega)$ and $g\in L^2(N)$ there exists a
\textit{unique} solution $u$ to the mixed problem
\eqref{mixedproblem}.  In addition there is a constant $C$
independent of $u$ such that
$$\int_{\partial \Omega} (\nabla u^*)^2 ds \leq C\left(\|f\|^2_{D^o}+\int _N g^2 ds\right)$$
\end{thm}

In the following section it is proved that a change from Dirichlet
to Neumann data on a single face is necessarily prohibited when the
change takes place across the graph of a Lipschitz function.  The
strict convexity condition of \eqref{partition} is also shown to be
necessary.  In the final section compact polyhedra are discussed for
which the set $N$ is necessarily \textit{empty}.

\section{Proof of Theorem \ref{mixed}}\mylabel{proof}

The estimates that follow are scale invariant.  Therefore to lighten
the exposition a bit it will be assumed, when working near any
vertex of the boundary of the compact polyhedron
$\overline{\Omega}$, that the vertex is at least a distance of $4$
units from any other vertex. Because $\partial \Omega$ is assumed to
be a $2$-manifold it will also be assumed that each edge that does
not contain a given vertex $v$ as an endpoint is at least $4$ units
from $v$ and similarly each face. Consequently, by another
application of the manifold condition, the picture that emerges is
that the \textit{truncated cones}

$$\mathcal{C}(v,r)=\{X\in\overline{\Omega}: |v-X|\leq r\}$$
for any vertex $v$ and $0\leq r\leq 4$ are homeomorphic to the
closed ball $\mathbb{B}^3$ while the \textit{cone bases}
$$\mathcal{B}(v,r)=\{X\in\overline{\Omega}: |v-X|= r\}$$
are homeomorphic to the closed disc $\mathbb{B}^2$.

Define
$$\Omega_r=\Omega\setminus \bigcup_{v}\mathcal{C}(v,r),\; 0<r<2$$
where the finite union is over all boundary vertices.  Then each
$\Omega_r$ is a Lipschitz domain (see,  for example, \S12.1 of
\cite{VV06} and Theorem 6.1 of \cite{VV03} for a proof and
generalizations in dimensions $n\geq3$). Likewise the interiors of
the \textit{arches} defined by
\begin{equation}\mylabel{arch}\mathcal{A}(v,r,R)=\{X\in\overline{\Omega}:r\leq |v-X|\leq
R\}\mbox{ ,    } 0<r<R<4\end{equation} are Lipschitz domains.

In general neither of these kinds of domains have a uniform
Lipschitz nature as $r\rightarrow 0$.  Therefore the following
\textit{polyhedral Rellich identity} of \cite{VV06} will be of use.
It is proved as in \cite{JK81} by an application of the Gauss
divergence theorem, but with respect to the vector field
$$W:= \frac{X}{|X|}, X\in \mathbb{R}^3\setminus \{0\}$$ when the
origin is on the boundary of the domain.
\begin{lem}\mylabel{rellichlemma}
Let $\mathcal{A}$ be any \emph{arch} \eqref{arch} of the polyhedral
domain $\Omega\subset \mathbb{R}^3$ and suppose $u$ is harmonic in
$\mathcal{A}$ with $\nabla u ^*\in L^2(\partial  \mathcal{A})$.
Then, taking the vertex $v$ to be at the \emph{origin}
\begin{equation}\mylabel{rellich}
2\int_ \mathcal{A} (W\cdot\nabla u)^2 \frac{dX}{|X|}=\int_
{\partial\mathcal{A}} \nu\cdot W|\nabla u|^2-2\partial_\nu u W\cdot
\nabla u ds
\end{equation}
\end{lem}

\begin{lem}\mylabel{rellichestimate}
With $\mathcal{A}=\mathcal{A}(v,r,R)$  and $u$ as in Lemma
\ref{rellichlemma}
\begin{equation}
2\int_ \mathcal{A} (W\cdot\nabla u)^2 \frac{dX}{|X|}\leq \int_
{\mathcal{B}(v,R)} |\nabla u|^2ds + 2\int_ {\mathcal{B}(v,r)}
(W\cdot\nabla u)^2ds+2\int_{\partial\Omega \cap \mathcal{A}}
|\partial_\nu u| |\nabla_t u| ds
\end{equation}
\end{lem}

\begin{proof}
The term $\nu\cdot W$ on the right of \eqref{rellich} is negative on
$\mathcal{B}(v,r)$ and vanishes on $\partial \Omega$.  Likewise the
second integrand on the right of \eqref{rellich} is a perfect square
on $\mathcal{B}(v,r)$, the negative of a square on
$\mathcal{B}(v,R)$, and $W\cdot \nabla u$ is a tangential derivative
on $\partial\Omega$.
\end{proof}

The partition $D\cup N=\partial\Omega$ induces a decomposition of
the Lipschitz boundaries $\partial\Omega_r$ into a Dirichlet part, a
Neumann part, and bases $\mathcal{B}(v,r)$ of the cones removed from
$\Omega$.  Define

$$N_r=(N\cap\partial\Omega_r)\bigcup_v \mathcal{B}(v,r)$$ and
$$D_r=\partial\Omega_r\setminus N_r.$$  \textit{This partition of
$\partial\Omega_r$ satisfies the requirements of a creased domain
in} \cite{Bro94}.  See \cite{VV06} pp. 586-587.  (Including the
bases in the Dirichlet part would also satisfy the requirements.) It
will therefore be possible to invoke Brown's \textit{existence}
results in the domains $\Omega_r$.

Similarly, \textit{arches $\mathcal{A}=\mathcal{A}(v,r,R)$ are
creased Lipschitz domains} with
$$N_r^R(v)=(N\cap\partial\mathcal{A}(v,r,R))\cup\mathcal{B}(v,r)\cup\mathcal{B}(v,R)$$
and
$$D_r^R=\partial\mathcal{A}\setminus N_r^R$$ for each vertex $v$.

Brown's estimate from \cite{Bro94} Theorem 2.1 is not scale
invariant.  However, the following special case is.

\begin{thm}\mylabel{brownnodirichlet}
\emph{(R. M. Brown)} Let $G\subset\mathbb{R}^n$ be a \emph{creased}
Lipschitz domain with $\partial G=D\cup N$. Then there exists a
unique solution $u$ to the mixed problem \eqref{mixedproblem} for
data $f$ identically \emph{zero} and $g\in L^2(N)$.  Furthermore
there is a constant $C$ determined only by the scale invariant
geometry of $G, D \mbox{ and } N$ and independent of $g$ such that
$$\int_{\partial G} (\nabla u^*)^2 ds \leq C\int _N g^2 ds$$
\end{thm}
\noindent As is
\begin{thm}\mylabel{brownconnected}
\emph{(R. M. Brown)} Let $G\subset\mathbb{R}^n$ be a \emph{creased}
Lipschitz domain with $\partial G=D\cup N$. Suppose that $D$ is
\emph{connected}.  Then there is a constant $C$ such that for all
harmonic functions  $u$ with $\nabla u^*\in L^2(\partial G)$
$$\int_{\partial G} (\nabla u^*)^2 ds \leq C\left(\int _D |\nabla_t u|^2 ds+\int _N (\partial_\nu u)^2 ds\right)$$
\end{thm}

\begin{proof}
Subtracting from $u$ its mean value  over $D$  allows the Poincar\'e
inequality over the connected set $D$. The conclusion still applies
to $u$.
\end{proof}

\begin{lem}\mylabel{nodirichlet}
Let $\Omega\subset \mathbb{R}^3$ be a compact polyhedral domain with
$2$-manifold boundary partitioned as $\partial\Omega= D\cup N$.  Let
$v$ be a vertex and let $j$ be a natural number.  Suppose $u$ is
harmonic in the arch $\mathcal{A}(v,2^{-j},2)$ with $\nabla u^*\in
L^2(\partial \mathcal{A})$ and $u$ \emph{vanishing} on
$D^2_{2^{-j}}$. Then there is a constant $C$ independent of $j$ so
that
$$\int_{\partial \Omega\cap\mathcal{A}(v,2^{-j},2) } |\nabla u|^2 ds \leq C\left(\int _{\partial\Omega\cap N_{2^{-j}}} (\partial_\nu u)^2 ds+
\int_ {\mathcal{B}(v,2^{-j})} (W\cdot\nabla u)^2ds +
\int_{\mathcal{A}(v,1,2) } |\nabla u|^2 dX\right)$$
\end{lem}

\begin{proof}
For natural numbers $k\leq j$ and real numbers $1\leq t\leq2$ the
arches $\mathcal{A}_{k,t}:=\mathcal{A}(v,t2^{-k},t2^{1-k})$ are
geometrically \textit{similar} Lipschitz domains.  Therefore by the
scale invariance  of Brown's  Theorem \ref{brownnodirichlet} above
$$\int_{\partial \Omega\cap\mathcal{A}_{k,t}} |\nabla u|^2 ds \leq C\int _{ N_{t2^{-k}}^{t2^{1-k}}} (\partial_\nu u)^2
ds$$ with $C$ independent of $k$.  Take $v$ to be the origin.  For
each $k$, integrating in $1\leq t\leq2$ and observing that $\nu= W$
or $-W$ on any cone base $\mathcal{B}$
$$\frac{1}{2}\int_{\partial \Omega\cap(\mathcal{A}_{k,1}\cup \mathcal{A}_{k,2})} |\nabla u|^2 ds
 \leq 2C\left(\int _{ N\cap(\mathcal{A}_{k,1}\cup \mathcal{A}_{k,2})} (\partial_\nu u)^2
ds+\int_ {\mathcal{A}_{k,1}\cup \mathcal{A}_{k,2}} (W\cdot\nabla
u)^2 \frac{dX}{|X|}\right)$$

Summing on $k=1,2,\ldots,j$ and using Lemma \ref{rellichestimate} on
the arch $\mathcal{A}(v,2^{-j},R)$ for each $1\leq R \leq 2$
together with the vanishing of $u$ on $D_{2^{-j}}^2$ again
\begin{multline*}\frac{1}{2}\int_{\partial \Omega\cap\mathcal{A}(v,2^{-j},2) } |\nabla u|^2 ds \leq
4C(\int _ {\partial\Omega\cap N_{2^{-j}}} (\partial_\nu u)^2 ds
+\\\int_ {\mathcal{B}(v,R)} |\nabla u|^2ds + 2\int_
{\mathcal{B}(v,2^{-j})} (W\cdot\nabla u)^2ds+2\int_{\partial\Omega
\cap N_{2^{-j}}^R} |\partial_\nu u| |\nabla_t u| ds+
\int_{\mathcal{A}(v,1,2) } |\nabla u|^2 dX)\end{multline*} An
application of Young's inequality
($2ab\leq\frac{1}{\epsilon}a^2+\epsilon b^2$) allows the square of
the tangential derivatives in the second to last term to be hidden
on the left side and the normal derivatives to be incorporated in
the first right side integral. Integrating in $1\leq R\leq2$ yields
the final inequality.
\end{proof}

By the same arguments, but using Theorem \ref{brownconnected} and
then Young's inequality in suitable ways for the $D$ portion and the
$N$ portion of the last integral of Lemma \ref{rellichestimate}, the
next lemma is proved.  For a given vertex, $D\cap
\mathcal{C}(v,R_1)$ is connected if and only if any $D\cap
\mathcal{A}(v,r,R_2)$ is connected.

\begin{lem}\mylabel{connected}
Let $\Omega\subset \mathbb{R}^3$ be a compact polyhedral domain with
$2$-manifold boundary partitioned as $\partial\Omega= D\cup N$.  Let
$v$ be a vertex and let $j$ be a natural number.  Suppose
$D\cap\mathcal{C}(v,2)$ is \emph{connected} and $u$ is harmonic in
the arch $\mathcal{A}(v,2^{-j},2)$ with $\nabla u^*\in L^2(\partial
\mathcal{A})$. Then there is a constant $C$ independent of $j$ so
that \begin{multline*}\int_{\partial
\Omega\cap\mathcal{A}(v,2^{-j},2) } |\nabla u|^2 ds \leq\\
C\left(\int _{D_{2^{-j}}} |\nabla_t u|^2 ds+\int
_{\partial\Omega\cap N_{2^{-j}}} |\partial_\nu u|^2 ds+ \int_
{\mathcal{B}(v,2^{-j})} (W\cdot\nabla u)^2ds +
\int_{\mathcal{A}(v,1,2) } |\nabla u|^2 dX\right)\end{multline*}
\end{lem}

Let $v$ be a vertex of the compact polyhedral domain $\Omega$ and
consider the collection of nontangential approach regions
$\Gamma(P)$ for $G=\Omega$ and parameter $\alpha$ \eqref{approach}
with $P\in\partial\Omega\cap\mathcal{C}(v,4)$.  By scale invariance
each approach region can be \textit{truncated} to a region
 $$\Gamma^T(P)=\{X\in \Gamma(P):|X-P|<(1+\alpha) dist(X,\partial\mathcal{A}(v,r/2,2r) )\},\;\;|v-P|=r$$
 so that the collections $\{\Gamma^T(P):r\leq |v-P|\leq 2r\}$ can be
 extended in a uniform way to systems of nontangential approach
 regions \textit{regular} in the sense of Dahlberg \cite{Dah79} for the
 arches $\mathcal{A}(v,r/2,4r)$.

 \textit{Denote by $w^T$ the nontangential maximal function of $w$ with respect to the truncated cones $\Gamma^T$.}

 \textit{Denote the Hardy-Littlewood maximal operator on $\partial \Omega$ by
 $\mathcal{M}$.} See, for example, \cite{Ste70} pp.10-11 or \cite{VV03}
 pp.501-502 for polyhedra.

 For $\alpha$ large enough a geometric argument shows that there is
 a constant independent of $P$ and $w$ such that
\begin{equation}\mylabel{hardylittlewood}w^*(P)\leq C \mathcal{M}(w^{T})(P)+ \max_K
|w|,\;\;P\in\partial\Omega\cap\mathcal{C}(v,4)\end{equation} where
$K$ is a compactly contained set in the Lipschitz domain $\Omega_2$.

Using Theorems \ref{brownnodirichlet} and \ref{brownconnected} to
estimate the truncated maximal functions introduces into the proofs
of Lemmas \ref{nodirichlet} and \ref{connected} a doubling of the
dyadic arches and therefore one dyadic term that is not immediately
hidden by Young's inequality.  Thus by the same proofs

\begin{lem}\mylabel{nodirichlet2}
With the same hypotheses as Lemma \ref{nodirichlet} there is a
constant $C$ independent of $j$ so that
\begin{multline*}\int_{\partial \Omega\cap\mathcal{A}(v,2^{1-j},2) } \left(\nabla u^T\right)^2 ds
-\frac{1}{2}\int_{\partial \Omega\cap\mathcal{A}(v,2^{-j},2^{1-j}) }
|\nabla u|^2 ds \leq \\C\left(\int _{\partial\Omega\cap N_{2^{-j}}}
(\partial_\nu u)^2 ds+ \int_ {\mathcal{B}(v,2^{-j})} (W\cdot\nabla
u)^2ds + \int_{\Omega_1 } |\nabla u|^2 dX\right)\end{multline*}
\end{lem}

\begin{lem}\mylabel{connected2}
With the same hypotheses as Lemma \ref{connected} there is a
constant $C$ independent of $j$ so that
\begin{multline*}\int_{\partial \Omega\cap\mathcal{A}(v,2^{1-j},2) }
\left(\nabla u^T\right)^2 ds-\frac{1}{2}\int_{\partial
\Omega\cap\mathcal{A}(v,2^{-j},2^{1-j}) } |\nabla u|^2 ds \leq\\
C\left(\int _{D_{2^{-j}}} |\nabla_t u|^2 ds+\int
_{\partial\Omega\cap N_{2^{-j}}} |\partial_\nu u|^2 ds+ \int_
{\mathcal{B}(v,2^{-j})} (W\cdot\nabla u)^2ds + \int_{\Omega_1 }
|\nabla u|^2 dX\right)\end{multline*}
\end{lem}

\begin{rem}\mylabel{lemmaremark}
Lemmas \ref{nodirichlet} and \ref{connected} apply to the negative
terms of Lemmas \ref{nodirichlet2} and \ref{connected2}.
Consequently those terms may be removed from the inequalities.
\end{rem}

\subsection{The regularity problem}

The \textit{regularity problem} is the mixed problem for
$\partial\Omega=D$.

\begin{thm}\mylabel{regularity}
Let $\Omega\subset\mathbb{R}^3$ be a compact polyhedral domain with
$2$-manifold connected boundary.  Then for any $f\in
W^{1,2}(\partial\Omega)$ the regularity problem is uniquely solvable
and the estimate for the solution $u$ $$\int_{\partial \Omega}
|\nabla u^*|^2 ds \leq C\int _{\partial \Omega} |\nabla_t f|^2 ds$$
holds with $C$ independent of $f$.
\end{thm}

\begin{proof}
For each $\Omega_{2^{-j}}$ there is unique solution $u_j$ to the
mixed problem with $u_j=f$ on $D_{2^{-j}}$ and $\partial_\nu u_j=0$
on $N_{2^{-j}}$ by Brown's existence result \cite{Bro94}.  By
definition of the truncated approach regions in each vertex cone
$\mathcal{C}(v,4)$ the regions may be extended to a regular system
of truncated approach regions for the $\partial \Omega\cap\partial
\Omega_1$ part of the boundary.  Thus the truncated nontangential
maximal function can be defined there. By Lemma \ref{connected2} and
Remark \ref{lemmaremark}, summing over all vertices, using analogous
estimates on the local Lipschitz boundary of $\partial\Omega$
outside of the vertex cones and using $W\cdot\nabla u_j=0$ on the
bases $\mathcal{B}(v,2^{-j})$,
\begin{equation}\mylabel{useconnected2}\int_{D_{2^{1-j}}}
\left(\nabla u^T_j\right)^2 ds \leq\\
C\left(\int _{D_{2^{-j}}} |\nabla_t f|^2 ds+ \int_{\Omega_1 }
|\nabla u_j|^2 dX\right)\end{equation} with $C$ independent of $j$.

Subtracting from $u_j$ the mean value $m_f$ of $f$ over
$\partial\Omega$ does not change \eqref{useconnected2}. Thus
Poincar\'e (see \cite{VV06}p.639 for polyhedral boundaries) can be
applied over $\partial\Omega$ with constant independent of $j$ in
\begin{equation}\mylabel{poincare}\int_{\Omega_1 }
|\nabla u_j|^2 dX\leq \int_{D_{2^{-j}}}  (u_j-m_f)\partial_\nu u_j ds=\int_{D_{2^{-j}}}  (f-m_f)\partial_\nu u_j ds \leq\\
C_\epsilon\int _{\partial \Omega} |\nabla_t f|^2 ds+\epsilon
\int_{D_{2^{-j}}} |\nabla u_j|^2 ds\end{equation} Applying Lemma
\ref{connected} to the part of the integral over the regions
$D_{2^{-j}}^{2^{1-j}}$ and using $W\cdot\nabla u_j=0$ again
$$\epsilon
\int_{D_{2^{-j}}} |\nabla u_j|^2 ds\leq \epsilon C\left(\int
_{D_{2^{-j}}} |\nabla_t f|^2 ds+ \int_{ \Omega_1} |\nabla u_j|^2
dX\right)+\epsilon \int_{D_{2^{1-j}}} |\nabla u_j|^2 ds$$ so that
\eqref{poincare} yields $$\frac{1}{2}\int_{\Omega_1 } |\nabla u_j|^2
dX\leq (C_\epsilon+\epsilon C)\int _{\partial \Omega} |\nabla_t f|^2
ds+\epsilon \int_{D_{2^{1-j}}} |\nabla u_j|^2 ds$$ for all
$\epsilon$ chosen small enough depending on $C$ but not on $j$.
Using this in \eqref{useconnected2} for $\epsilon$ chosen small
enough gives
\begin{equation}\mylabel{regest}
\int_{D_{2^{1-j}}}
\left(\nabla u^T_j\right)^2 ds \leq\\
C'\int _{\partial\Omega} |\nabla_t f|^2 ds
\end{equation}
with the constant independent of $j$.

Given any compact subset of $\Omega$, \eqref{regest} together with
$u_j=f$ on $D_{2^{-j}}$ for all $j$ implies there exists a
subsequence so that both $u_{j_k}$ and $\nabla u_{j_k}$ converge
uniformly on the compact set to a harmonic function $u$ and its
gradient respectively. A diagonalization argument gives pointwise
convergence on all of $\Omega$.  Intersecting a compact  subset $K$
with the truncated approach regions yields compactly contained
regions and corresponding maximal functions $\nabla
u_{j_k}^{T,K}\rightarrow\nabla u^{T,K} $ uniformly. Thus by
\eqref{regest} and then monotone convergence, as $\Omega$ is
exhausted by compact subsets $K$,
\begin{equation}\mylabel{regest2}
\int_{\partial\Omega}
\left(\nabla u^T\right)^2 ds \leq\\
C\int _{\partial\Omega} |\nabla_t f|^2 ds
\end{equation}
See \cite{JK82b} for these arguments.

A difficulty with the setup here is that the $\int_{D_{2^{1-j}}}
\left(\nabla(u_j- u_k)^T\right)^2 ds$ for $k>j$ do not \textit{a
priori} have better bounds than the right side of \eqref{regest}.
However, \eqref{regest} together with weak convergence in
$L^2(\partial\Omega_{2^{-j}})$ and pointwise convergence on the
bases $\mathcal{B}(v,2^{-j})$ shows that for each $j$ and every
$X\in\Omega_{2^{-j}}$    a subsequence of
$$u_k(X)=\int_{\partial\Omega_{2^{-j}}} \mathcal{P}^X_j u_k ds=\int_{D_{2^{-j}}}
\mathcal{P}^X_j f ds+\sum_v\int_{\mathcal{B}(v,2^{-j})}
\mathcal{P}^X_j u_k ds$$ converges to $u(X)$, perforce with Poisson
representation that must be an extension from $D_{2^{-j}}$ of $f$.
Here $ \mathcal{P}^X_j$ is the Poisson kernel for the Lipschitz
polyhedral domain $\Omega_{2^{-j}}$ and may be seen to be in
$L^2(\partial\Omega_{2^{-j}})$ by Dahlberg \cite{Dah77}.
Consequently $u$ has nontangential limits $f$ on $\partial\Omega$,
and by \eqref{hardylittlewood} and \eqref{regest2} the theorem is
proved.
\end{proof}

\subsection{The mixed problem with vanishing Dirichlet data}

\begin{thm}\mylabel{v}
Let $\Omega\subset\mathbb{R}^3$ be a compact polyhedral domain with
$2$-manifold connected boundary.  Then for any $g\in L^2(N)$ there
is a unique solution $u$ to the mixed problem \eqref{mixedproblem}
that vanishes on $D$ and has Neumann data $g$ on $N$.  Further
$$\int_{\partial \Omega} (\nabla u^*)^2 ds \leq C\int _N g^2 ds$$
\end{thm}

\begin{proof}
Again by \cite{Bro94} there exists a unique solution $u_j$ in
$\Omega_{2^{-j}}$ to the mixed problem so that $\partial_\nu u_j=g$
on $\partial \Omega\cap N_{2^{-j}}$, $W\cdot\nabla u_j=0$ on the
$\mathcal{B}(v,2^{-j})$ and $u_j=0$ on $D_{2^{-j}}$.  Lemma
\ref{nodirichlet2} and Remark \ref{lemmaremark} imply
$$\int_{\partial \Omega\cap\partial \Omega_{2^{1-j}} } \left(\nabla u_j^T\right)^2 ds
 \leq \\C\left(\int _{\partial\Omega\cap N_{2^{-j}}}
g^2 ds+  \int_{\Omega_1 } |\nabla u_j|^2 dX\right)$$

A Poincar\'e inequality independent of $j$ is also needed here and
is supplied by the following lemma. Polyhedral domains are naturally
described as simplicial complexes. See for definitions and notations
\cite{Gla70} \cite{RS72} \cite{VV03}
 \cite{VV06} or others.
\begin{lem}  Suppose $u$ is harmonic in $\Omega_{2^{-j}}$ with
$\partial_\nu u=g$ on $\partial \Omega\cap N_{2^{-j}}$,
$\partial_\nu u=0$ on the $\mathcal{B}(v,2^{-j})$ and $u=0$ on
$D_{2^{-j}}$.  Then
$$\int_{ \Omega_{2^{-j}} } \left|\nabla u\right|^2 dX
 \leq \\C\int _{\partial\Omega\cap N_{2^{-j}}}
g^2 ds$$ with $C$ independent of $j$.
\end{lem}

\begin{proof}  By Green's first identity and Young's inequality
\begin{equation}\mylabel{poincare2}\int_{\Omega_{2^{-j}} }
|\nabla u|^2 dX= \int_{\partial\Omega\cap N_{2^{-j}}}  u\;\partial_\nu u ds \leq\\
C_\epsilon\int _{\partial\Omega\cap N_{2^{-j}}} g^2 ds+\epsilon
\int_{\partial\Omega\cap N_{2^{-j}}} u^2 ds\end{equation}

The polyhedron $\overline{\Omega}$ can be realized as a finite
homogeneous simplicial $3$-complex.  A cone $\mathcal{C}(v,1)$ is
then the intersection of the ball $|X|\leq 1$ with the \textit{star}
$St(v,\overline{\Omega})$ in the $3$-complex $\overline{\Omega}$ of
the vertex $v$. Each $2$-simplex $\sigma^2$ of
$St(v,\overline{\Omega})$ that is also contained in $N$ is contained
in a unique $3$-simplex $\sigma^3\in St(v,\overline{\Omega})$.  Let
$B$ denote the \textit{unit} vector in the direction from the
\textit{barycenter} of $\sigma^3$ to $v$.  Then $\sigma^2\cap
\{|X|\leq 1\}$ may be projected into the sphere $|X|=1$ along lines
parallel to $B$ by $Q\mapsto Q+t_Q B$ onto a set contained in
$\sigma^3\cap \mathcal{B}(v,1)$.  The sets $\{Q+tB: Q\in\sigma^2\cap
N^1_{2^{-j}}(v)\mbox{ and } 0\leq t\leq t_Q\}$  are contained in
$\sigma^3\cap\mathcal{A}(v,2^{-j},1)$.  Thus by the fundamental
theorem of calculus for each $Q\in\partial\Omega\cap
N^1_{2^{-j}}(v)$ and integrating $ds(Q)$
\begin{equation}\mylabel{conepoincare}
\int_{\partial\Omega\cap N^1_{2^{-j}}(v)}  u^2 ds \leq C\left(\int
_{\mathcal{A}(v,2^{-j},1)} |\nabla u|^2 dX+\int _{\mathcal{B}(v,1)}
u^2 ds\right)
\end{equation}
where the constant depends only on the projections, i.e. only on the
finite geometric properties of the complex that realizes
$\overline{\Omega}$ and not on $j$.

By the fundamental theorem, the connectedness of $\Omega_1$ and the
vanishing of $u$ on the fixed nonempty set $D$
$$\int_{\partial\Omega_1} u^2 ds \leq C\int _{\Omega_1} |\nabla u|^2
dX$$ This together with \eqref{conepoincare} implies
$$\epsilon
\int_{\partial\Omega\cap N_{2^{-j}}} u^2 ds\leq \epsilon C \int_{
\Omega_{2^{-j}} } \left|\nabla u\right|^2 dX$$ and $\epsilon$ can be
chosen independently of $j$ so that \eqref{poincare2} yields the
lemma.
\end{proof}

The lemma yields the analogue of \eqref{regest}
\begin{equation}\mylabel{regest2}
\int_{\partial \Omega\cap\partial \Omega_{2^{1-j}} } \left(\nabla
u_j^T\right)^2 ds
 \leq \\C\int _{\partial\Omega\cap N_{2^{-j}}}
g^2 ds
\end{equation}
Continuing to argue as in the proof of Theorem \ref{regularity} ,
this and the vanishing of the $u_j$ on $D_{2^{-j}}$ produces a
harmonic function $u$ defined in $\Omega$ that is the pointwise
limit of a subsequence of the $u_j$.  In addition $u$ satisfies
$$\int_{\partial \Omega } \left(\nabla
u^T\right)^2 ds
 \leq \\C\int _{ N}
g^2 ds$$ which in turn yields the maximal estimate of the theorem.

To show that $u$ assumes the correct data, \eqref{regest} along with
weak $L^2$-convergence, pointwise convergence and the Poisson
representation in each $\Omega_{2^{-j}}$ proves as before that $u$
vanishes nontangentially on $D$.  By constructing a Neumann function
(possible by \cite{JK81}) in analogy to the Green function, or by
using the invertibility of the classical layer potentials
\cite{Ver84}, a Neumann representation of $u$ in each
$\Omega_{2^{-j}}$ can be obtained so that $\frac{\partial
u}{\partial \nu}=g$ nontangentially on $N$ can be deduced by the
same arguments.

Uniqueness follows from Green's first identity valid in polyhedra
when $\nabla u^*\in L^2$.
\end{proof}

\subsection{Proof of Theorem \ref{mixed}}

Recall the definition of the homogeneous Sobolev semi-norm
\eqref{seminorm}.

\begin{lem}
When $\partial\Omega$ is connected and $\|f\|^2_{D^o}=0$,  $f$ is
identically constant on $D$.
\end{lem}

The lemma says that $f$ equals the same constant value on each
component of $D$.

\begin{proof}
Because the semi-norm equals \textit{zero} there is a sequence of
extensions $\widetilde{f}_j$ of $f$ and a sequence of numbers $m_j$
so that by Poincar\'e in the second inequality
$$\int_D(f-m_j)^2 ds\leq\int_{\partial\Omega}(\widetilde{f}_j-m_j)^2 ds\leq C\int _{\partial\Omega} |\nabla_t \widetilde{f}_j|^2 ds\rightarrow 0$$
\end{proof}

\noindent \textit{Proof of Theorem \ref{mixed}.}  Choose an
extension $\widetilde{f}$ of $f$ so that $\int _{\partial\Omega}
|\nabla_t \widetilde{f}|^2 ds \leq 2 \|f\|^2_{D^o}$.  This is always
possible by the lemma.  Then from Theorem \ref{regularity} there is
a unique solution $u_D$ with regularity data $\widetilde{f}$ and
$\int_{\partial \Omega} (\nabla u_D^*)^2 ds \leq C\|f\|^2_{D^o}$.
From Theorem \ref{v} there is a unique solution $u_N$ vanishing on
$D$, with Neumann data $g-\partial_\nu u_D$ on $N$, and
$$\int_{\partial \Omega} (\nabla u_N^*)^2 ds \leq
C\left(\int _N (\partial_\nu u_D)^2 ds+\int _N g^2 ds\right)\leq
C\left(\|f\|^2_{D^o}+\int _N g^2 ds\right)$$ The solution is
$u=u_D+u_N$.  Theorem \ref{v} established uniqueness. $\Box$

\section{On violations of the postulates for the partition
$\partial\Omega=D\cup N$}

When $D$ is empty the mixed problem is the Neumann problem and
solvable for any data that has mean value \textit{zero} on the
boundary \cite{Ver01}.  We consider the two remaining postulates.

\subsection{$N$ is the union of a number (possibly \textit{zero}) of
closed faces of $\partial\Omega$.}\mylabel{subsec}

Solving the mixed problem means that every $W^{1,2}(D)$ function has
a $W^{1,2}(\partial \Omega)$ \textit{extension}.  This observation
raises the possibility that the mixed problem might be solvable when
a given (open) face $\mathcal{F}$ has  nonempty intersection both
with $D$ and with $N$ in such a way that $D\cap \mathcal{F}$ is an
\textit{extension domain}.  Here we will only consider the
possibility that this extension domain has a Lipschitz boundary
\cite{Ste70} and show that \textit{the mixed problem is never
solvable when this condition on the partition occurs}.

 Let $\phi : \mathbb{R}\rightarrow
\mathbb{R}$ be a Lipschitz continuous function $y=\phi(x)$ with
$\|\phi'\|_\infty \leq M$. Choose a point $p_0$ on the graph
$(x,\phi (x))$ in the plane and consider the rectangle with width
$2$ parallel to the $x$-axis, length $8M$ and \textit{center} $p_0$.
Locate the origin directly below $p_0$ and $M$ units from the bottom
of the rectangle.  Here it will be convenient to name the region $N$
that is \textit{strictly below} the graph and contained in the
rectangle.  Call its complement in the rectangle $D$.  Let $(x,y,t)$
be the rectangular coordinates of $\mathbb{R}^3$ with
\textit{origin} coinciding with the origin of the plane.  Let $Z$ be
the open right circular cylinder of $\mathbb{R}^3$ with center $p_0$
that intersects the plane in precisely the (open) rectangle.

The domain $\Omega=Z\setminus D \subset \mathbb{R}^3$ is
\textit{regular} for the Dirichlet problem.  This follows by the
Wiener test applied to each of the points of
$\partial\Omega=\partial Z \cup D$ and the observation that the
Newtonian capacity in $\mathbb{R}^3$ of a disc from the plane is
proportional to its radius.  See, for example, \cite{Lan72} p. 165.
Here the Lipschitz (or NTA) condition is also used.  Consequently
the Green function, $g=g^0$ for $\Omega$ with pole at the origin, is
\textit{continuous} in $\overline{\Omega}\setminus \{0\}$.

Approximating Lipschitz domains to $\Omega$ are constructed as
follows.  For each $\tau>0$ define Lipschitz surfaces with boundary
(the graph of $\phi$) by
$$D_\tau = \{p+s(p-\tau e_3): p \mbox{ is on the graph of } \phi
\mbox{ and } 0\leq s \} \cap Z$$ Here $e_3$ is the standard basis
vector perpendicular to the $xy$-plane.  Denote by $H_\tau$ the
region of $Z$ between $D$ and $D_\tau$ above the graph of $\phi$.
Then the $\Omega_\tau =\Omega\setminus \overline{H}_\tau=Z\setminus
\overline{H}_\tau$ are Lipschitz domains.  Denote by $g_\tau$ the
Green function for $\Omega_\tau$ with pole at $\tau e_3$.

\begin{lem}\mylabel{main}

(i) $-\partial_tg(x,y,t)$ for $t>0$ has continuous boundary values
$\partial_\nu g:= -\lim_{t\downarrow 0} g(x,y,t)/t$  at every point
of $D$ for which $y>\phi(x)$.

(ii) $\int_D (\partial _\nu g)^2 ds= +\infty$.

(iii) $\partial_tg(x,y,0)=0$ at every point of $N\setminus \{0\}$.

(iv) $\partial_\nu g \in L^2(\partial Z)$.
\end{lem}

\begin{proof}
(i) follows by Schwarz reflection while (iii) follows by the
symmetry in $t$ of $\Omega$ and $g$.  The maximum principle shows
that the Green function for $Z$ dominates from below the Green
function for $\Omega$, $g_Z\leq g\leq 0$.  On $\partial Z$ both
Green functions vanish so that $\partial_\nu g_Z\geq \partial_\nu
g\geq 0$ while $\partial_\nu g_Z$ is square integrable there,
establishing (iv).

D. S. Jerison and C. E. Kenig's Rellich identity for harmonic
measure (\cite{JK82b} Lemma 3.3) is valid on any \textit{Lipschitz}
domain $G$ that contains the \textit{origin}.  It is
$$(n-2)w_G(0)=\int_{\partial G}(\partial_\nu g_G (Q))^2 \nu\cdot Q
ds(Q)$$ with respect to the vector field $X$.  Here
$g_G(X)=F(X)+w_G(X)$ is the Green function for $G$, and $F$ is the
fundamental solution for Laplace's equation.  Denote by $w_\tau$,
$w$ and $w_Z$ the corresponding harmonic functions for the
$\Omega_\tau$, $\Omega$ and $Z$ Green functions respectively.

By $Z\supset Z\setminus D=\Omega\supset \Omega_\tau$ and the maximum
principle
$$\partial_\nu g_\tau\leq \partial_\nu
g_Z \mbox{ on } \partial\Omega_\tau\setminus D\setminus D_\tau$$
$$0< \partial_\nu g_\tau \leq g \mbox{ on } D$$
and
\begin{equation}\mylabel{w}
w_Z<w<w_\tau \mbox{ in } \Omega_\tau
\end{equation}
For $Q\in D$ and $\nu=\nu_Q$ the outer unit normal to $\Omega_\tau$,
$\nu\cdot(Q-\tau e_3)=\tau$, while for $Q\in D_\tau$,
$\nu\cdot(Q-\tau e_3)=0$.  Formulating the Rellich identity with
respect to the vector field $X-\tau e_3$ and using these facts
($n=3$)
\begin{multline*}
w_\tau(\tau e_3) = \int_{ \partial\Omega_\tau \setminus D \setminus
D_\tau}(\partial_\nu g_\tau )^2 \nu\cdot (Q-\tau e_3)
ds+\tau\int_{D}(\partial_\nu g_\tau )^2  ds\leq \\\int_{
\partial Z}(\partial_\nu
g_Z )^2 \nu\cdot (Q-\tau e_3) ds+\tau\int_{D}(\partial_\nu g )^2
ds=w_Z(\tau e_3) +\tau\int_{D}(\partial_\nu g )^2 ds
\end{multline*}
so that
$$\frac{w(\tau e_3)- w_Z(\tau e_3)}{\tau} <\frac{w_\tau(\tau e_3)- w_Z(\tau e_3)}{\tau} \leq\int_{D}(\partial_\nu g )^2
ds$$ and (ii) follows from \eqref{w} and $\tau \downarrow 0$.
\end{proof}

For $\delta >0$ define smooth subdomains of $\Omega$
$$G_\delta=\{g<-\delta \}.$$
$\partial G_\delta\rightarrow \partial\Omega$ uniformly.  The
$\partial_\nu g\mid_{ \partial G_\delta} ds$ are a collection of
probability measures on $\mathbb{R}^3$ that  have \textit{harmonic
measure for $\Omega$ at the origin} as  weak-$*$ limit.

By $G_\delta\uparrow \Omega$, Green's first identity, and monotone
convergence
\begin{equation}\mylabel{l2omega}
\int_{\Omega\setminus B_r}|\nabla g|^2 dX < \infty
\end{equation}
for all balls centered at the origin.

With $\phi$, $N$, $D$ and $Z$ as above define the half-cylinder
domain $Z_+=\{(x,y,t)\in Z: t>0\}$.  Then $D\cup N\subset \partial
Z_+\cap \{t=0\}$.

\begin{lem}\mylabel{Y}
Suppose $\triangle u=0$ in $Z_+$, $\nabla u^*\in L^2(\partial Z_+)$,
$\partial_\nu u\rightarrow^{n.t.} 0\; a.e.$ on $N$, and $
u\rightarrow^{n.t.} 0\; a.e.$ on $D$.  Let $Y\subset Z$ be a scaled
cylinder centered at $p_0$ with $dist(\partial Y, \partial Z)>0$.
Let $Y_+$ be the corresponding half-cylinder.  Then $u\in
C(\overline{Y}_+)$.

\end{lem}

\begin{proof}
The hypothesis on $\nabla u^*$ implies $u^*\in L^2(\partial Z_+)$ so
that $u$ and $\nabla u$ have nontangential limits $a.e.$ on
$\partial Z_+$ \cite{Car62} \cite{HW68}.  Extend $u$ to the bottom
component of $Z\setminus D\setminus N$ by $u(x,y,t)=u(x,y,-t)$.  By
the vanishing of the Neumann data on $N$, $\triangle u=0$ in the
sense of distributions in the domain $\Omega= Z\setminus D$ and then
classically.

Fix $d>0$ and suppose $X\in \overline{Y}_+$ is of the form
$X=(x,y,d)$ for $ y\geq \phi(x)-Md$.  Denote $3$-balls of radius and
 distance to $D$ \textit{comparable to}  $d$   by $B_d$.  Denote
 $2$-discs in $\partial Z_+$ with radius\textit{ comparable to} $d$ by
 $\Delta_d$ and let $\inta{}$ denote integral average.  Then by the
 mean value theorem, the fundamental theorem of calculus, the $a.e.$
 vanishing of $u$ on $D$ and the geometry of the nontangential
 approach regions
 $$ |u(X)|\leq \inta{B_d(X)}|u|\leq
 Cd(\inta{\Delta_d(x,y+2Md,0)}\nabla u^* ds)$$
 where $C$ depends only on $M$.  By absolute continuity of the
 surface integrals and $\nabla u^*\in L^2$ there is a function $\eta(d)\rightarrow 0$ as $
 d\rightarrow 0$ so that $\int_{\Delta_d} (\nabla u^{*})^2 ds\leq
 \eta(d)$ for all $\Delta_d\subset \partial Z_+$.  Consequently the
 Schwarz inequality now yields $|u(X)|\leq C \eta(d)$.

 Suppose now $X$ is of the form $X=(x, \phi(x)-Md,t)$ for $0\leq t
 \leq d$.  Because $u$ has been extended
 $$|u(X)|\leq \inta{B_d(X)}|u|\leq(\inta{B_d(x,\phi(x)-Md,d)}|u|)
 +d(\inta{\Delta_d(x,\phi(x)-Md,0)}\nabla u^* ds)$$
 and $|u(X)|\leq 2 C \eta(d)$.  The lemma follows.
\end{proof}

Partition $\partial Z_+$ by $N_+=\overline{N}$, $D_+=\partial
Z_+\setminus \overline{N}$ and $\partial Z_+=N_+\cup D_+$.  For
$3/4> r>0$ let $Z^r$ be the scaled cylinder centered at $p_0$ of
width $2r$ and length $8Mr$.  Define the corresponding
half-cylinders $Z_+^r$ with $$N_+^r=N_+\cap \partial Z^r_+$$ (not a
scaling of $N_+$) and $$D_+^r= \partial Z_+^r\setminus N_+^r$$ With
this partition  $Z_+^r$ is called a \textit{split cylinder with
Lipschitz crease}.

\textit{By \eqref{l2omega} and the Fubini theorem}, $g\in W^{1,2}
(\partial Z_+^r\setminus \{t=0\})$ \textit{for} $a.e.$ $r$.

\begin{prop}\mylabel{prop}
Let $g$ be the Green function for $\Omega=Z\setminus D$ with pole at
the origin.  For \textit{almost every} $\frac{3}{4}>r>0$ there
exists no solution $u$ to the $L^2$-\textit{mixed problem}
\eqref{mixedproblem} in the split cylinder with Lipschitz crease
$Z_+^r$ with boundary values $u\rightarrow^{n.t.} g\in
W^{1,2}(D_+^r)$ and $\partial_\nu u\rightarrow \partial_\nu g=0$ on
$N_+^r$.
\end{prop}

\begin{proof}
Suppose instead that there is such a solution $u$ with $\nabla
u^*\in L^2(\partial Z_+^r)$.  Then the first paragraph of the proof
of Lemma \ref{Y} applies and, in particular, $u$ extends to
$Z^r\setminus D$ evenly and harmonically across $N_+^r$.  The
Dirichlet data that $u$ takes $a.e.$ on $D_+^r$ is a continuous
function, as is the \textit{Dirichlet } data  that $u$ takes
(continuously) on $N_+^r$.  The Dirichlet data $u$ takes $a.e.$ on
$\partial Z_+^r$ will be shown to be a continuous function if it can
be shown to be continuous across the boundary $\partial N_+^r$ of
the surface $N_+^r$.  Lemma \ref{Y}, scaled to apply to the split
cylinders here, shows that the Dirichlet data is continuous across
the Lipschitz crease part of $\partial N_+^r$.  The same argument
used there works on the other parts:  Suppose $dist(X,\partial
Z^r)=d$ for $X\in N_+^r$.  Let $\Delta_d\subset \partial Z^r\cap
D_+^r$ be a disc approximately a distance $d$ from $X+de_3$.  Then
$$|u(X)-\inta{\Delta_d}g ds|\leq |\inta{B_d(X)}u(Y)-u(Y+de_3)dY| +
|\inta{B_d(X+de_3)}u(Y)-\inta{\Delta_d}g ds|\leq C\eta(d)$$ and the
continuity across $\partial N_+^r$ follows from the continuity of
$g$ and $\eta(d)\rightarrow 0$.

Thus the data $u$ takes $a.e.$ on $\partial Z_+^r$ is a continuous
function.  Since also $u^*\in L^2(\partial Z_+^r)$ it follows that
$u\in C(\overline{Z^r_+})$.  The evenly extended $u$ is then
continuous on $\overline{Z^r}$, harmonic in $Z^r\setminus D$ with
the same Dirichlet data as $g$ on $\partial (Z^r\setminus D)$.  The
maximum principle implies $u=g$.

Let $g_r$ denote the Green function for $Z^r\setminus D$ with pole
at a point $\{P\}$ of $N_+^r$.  Again $g_r$ is continuous in
$\overline{Z^r}\setminus \{P\}$.  Let $B\subset \overline{B}\subset
Z^r\setminus D$ be a ball centered at $P$.  Then by the maximum
principle $cg\geq g_r$ on $\overline{Z^r\setminus B}$ for some
constant $c$.  By this domination, the vanishing of both $g$ and
$g_r$ on $D_+^r\cap \{t=0\}$ and (ii) of Lemma \ref{main} applied to
$g_r$, it follows that $\partial_\nu g$ which is not in $L^2(D)$ can
neither be square integrable over the smaller set $D_+^r\cap
\{t=0\}$.  Since $u=g$ this contradicts the assumption on the
nontangential maximal function of the gradient.
\end{proof}

The nonsolvability of the $L^2$-mixed problem in the split cylinders
can be extended to nonsolvability in any polyhedron that has a
Lipschitz graph crease on any face by a  \textit{globalization
argument}. Let $g$ and $r$ be as in the Proposition.  By using the
approximating domains $Z^r\cap G_\delta$ as $\delta\rightarrow 0$,
the Green's representation
$$g(X)=\int_{\partial Z^r} \partial_\nu F^X g- F^X\partial_\nu g ds-
\int_{D\cap \overline{Z^r}} F^X d\mu^0, X\in Z^r\setminus D$$ can be
justified where $\mu^0$ is harmonic measure for $\Omega=Z\setminus
D$ at the origin and $F$ is the fundamental solution for Laplace's
equation.  Let $\chi\in C_0^\infty(\mathbb{R}^3)$ be a cut-off
function that is supported in a ball contained in $Z^r$ centered at
$p_0$, and is  identically $1$ in a concentric ball $B^r$ with
smaller radius. Then define $$u(X)= -\int_{D\cap \overline{Z^r}}F^X
\chi d\mu^0$$ harmonic in $\mathbb{R}^3$ outside $supp (\chi) \cap
D$. Similarly $g(X)-u(X)$ is harmonic inside $B^r$.  Consequently
$\nabla u^*\notin L^2(B^r\cap D)$ by applying a scaled (ii) of Lemma
\ref{main}  to $g$ again. Also
\begin{equation}\mylabel{urep}
u(X)= -\int_{D\cap \overline{Z^r}}F^X(Q)\left(
\chi(Q)-\chi(X)\right) d\mu^0(Q)-\chi(X)\int_{\partial Z^r}
\partial_\nu F^X g- F^X\partial_\nu g ds+\chi(X) g(X)
\end{equation}
The last term has bounded Neumann data on $N$ and vanishing
Dirichlet data on $D$.  The Cauchy data of the middle term is smooth
and compactly supported on $D\cup N$.  For any $X\notin D$ the
\textit{gradient} of the first term is bounded by a constant,
depending on $\chi$, times $$\left|\int_{D\cap \overline{Z^r}}F^X(Q)
d\mu^0(Q)\right|\leq -F^X(0)+g^X(0)\leq \frac{1}{4\pi|X|}$$ Here
$g^X$ is the (negative) Green function for $\Omega=Z\setminus D$
with pole at $X$.  Thus the first term is Lipschitz continuous on
$D_+^r\cup N_+^r$.  Altogether $u$ has bounded Neumann data on $N$
and Lipschitz continuous data on $D$ while $\nabla u^*\notin L^2(
D)$. Finally $\nabla u\in L^2_{loc}( \mathbb{R}^3)$ by \eqref{urep}
since this is true for $\chi g$.

Thus whenever a split cylinder $Z_+^r$ can be contained in a
polyhedral domain so that $\partial Z_+^r\cap \{t=0\}$ is contained
in a face and so that the Lipschitz crease is part of the boundary
between the Dirichlet and Neumann parts of the polyhedral boundary,
then the harmonic function $u$ just constructed is defined in the
entire polyhedra domain.  Its properties suffice to compare it with
any solution $w$ in the class $\nabla w^* \in L^2$ by Green's first
identity $\int |\nabla u-\nabla w|^2 dX= \int (u-w)\partial_\nu
(u-w) ds$.  Regardless of the nature of the partition away from
$Z_+^r$, when $w$ has the same data as does $u$ it must be
concluded, as in Proposition \ref{prop}, that it is identical to
$u$.  This establishes

\begin{thm}
Let $\Omega \subset \mathbb{R}^3$ be a compact polyhedral domain
with partition $\partial \Omega= D\cup N$.  Let $\mathcal{F}$ be an
open face of $\partial \Omega$ such that $\mathcal{F} \cap D$ is a
Lipschitz domain of $\mathcal{F}$ with nonempty complement
$\mathcal{F} \cap N$.  Then there exist mixed data
\eqref{mixedproblem} for which there are no solutions $u$ in the
class $\nabla u^*\in L^2(\partial \Omega)$.
\end{thm}
\subsection{Whenever a face of $N$ and a face of $D$ share a
$1$-dimensional edge as boundary, the dihedral angle measured in
$\Omega$ between the two faces is \textit{less} than $\pi$.}

Continue to denote points of $\mathbb{R}^3$ by $X=(x,y,t)$.  Define
$D$ to be the upper half-plane of the $xy$-plane.  Introduce polar
coordinates $y=r \cos \theta$ and $t=r \sin \theta$, let $\pi\leq
\alpha<2\pi$ and define $N$ to be the half-plane $\theta=\alpha$.
The crease is now the $x$-axis.

 Define
 $$b(X)=r^{\frac{\pi}{2\alpha}}\sin(\frac{\pi}{2\alpha}\theta)$$ for
 $X$ above $D\cup N$.  These are  Brown's counterexample
 solutions for nonconvex plane sectors \cite{Bro94}.  The Dirichlet
 data vanishes on $D$ while the Neumann vanishes on $N$, and $\nabla
 b^*\notin L^2$.

 These solutions are globalized to a compact polyhedral domain with
 interior dihedral angle $\alpha$:

 Denote by $\Theta$ the intersection of a (large) ball centered at
 the origin and the domain above $D\cup N$.  Then $b(X)$ is
 represented in $\Theta$ by $$b(X)=\int_{\partial \Theta \setminus
 D}\partial _\nu F^X b ds-\int_{\partial \Theta \setminus
 N} F^X \partial _\nu b ds$$
 Let $\chi\in C_0^\infty(\mathbb{R}^3)$ be a cut-off function as
 before , but centered at the \textit{origin} on the crease.  Define
$$u(X)=\int_{N}\partial _\nu F^X\chi b ds-\int_{D} F^X \chi\partial _\nu b ds$$
As before, $u$ is harmonic everywhere outside $supp (\chi) \cap
(D\cup N)$ and $\nabla u^* \notin L^2\left(supp (\chi) \cap (D\cup
N)\right)$. Also
\begin{multline}\mylabel{brep}
u(X)= \int_{N\cap \overline{\Theta}}\partial_\nu F^X(Q)\left(
\chi(Q)-\chi(X)\right) b(Q) ds(Q)\\-\int_{D\cap
\overline{\Theta}}F^X(Q)\left( \chi(Q)-\chi(X)\right)\partial_\nu
b(Q) ds(Q)-\chi(X)\int_{\partial \Theta\setminus N\setminus D}
\partial_\nu F^X b- F^X\partial_\nu b ds+\chi(X) b(X)
\end{multline}
Again the boundary values around the support of $\chi$ are the
issue.  The last two terms are described just as the middle and last
after \eqref{urep}.  The \textit{gradient} of the second term is
\textit{bounded} because the integral over $D$ can be no worse than,
for example, $\int_0^1dx\int_0^1
\frac{1}{\sqrt{x^2+r^2}}\frac{dr}{r^\beta}<\infty$ for any $\beta<1$
(e.g. $\beta = 1-\frac{\pi}{2\alpha})$.

For a $\frac{\partial}{\partial X_j}$ derivative define tangential
derivatives (in $Q$) to any surface with unit normal $\nu$ by
$\partial^t_i=\nu_i\partial_j-\nu_j\partial_i$.  Then by the
harmonicity of $F$ away from $X$ and the divergence theorem in
$\Theta$, the $\frac{\partial}{\partial X_j}$ derivative of the
first integral equals the sum in $i$ of$$\int_{N\cap
\overline{\Theta}}\partial_i F^X\partial_i^t\left(( \chi-\chi(X))
b\right) ds$$ plus integrals over $\partial\Theta\setminus
N\setminus D$ ($b$ vanishes on $D$) that will all be
\textit{bounded} since $X$ is near the support of $\chi$.  When the
tangential derivative falls on $b$ the integral is bounded like the
second term of \eqref{brep}.  The remaining integral has boundary
values in every $L^p$ for $p<\infty$ by singular integral theory.
(In fact, it too is bounded by a closer analysis, thus making it
consistent with the example from Section \ref{subsec}.)

Finally $\nabla u\in L^2_{loc}(\mathbb{R}^3)$ by its now established
properties and the corresponding property for $b$.   The argument
using Green's first identity as at the end of Section \ref{subsec}
is justified and\vskip.1in

\noindent\textit{The solutions $u$ can now be placed in polyhedral
domains that have interior dihedral angles greater than or equal to
$\pi$ and provide mixed data for which no $L^2$-solution can exist.}

\section{Polyhedral domains that admit only the trivial mixed
problem}

Consider the $L^2$-mixed problem for  the unbounded domain exterior
to a compact polyhedron.  When the polyhedron is convex the
requirement of postulate (iii) of \eqref{partition} eliminates all
but the trivial partition from the class of well posed mixed
problems.  In this case we will say that the \textit{exterior}
problem is \textit{monochromatic}.

The mixed problem for a compact polyhedral domain can also be
monochromatic for the \textit{interior} problem.  An example is
provided by the regular compound polyhedron that is the union of $5$
equal regular tetrahedra with a common center, a picture of which
may be found as Number $6$ on Plate III between pp.48-49 of H. S. M.
Coxeter's book \cite{Cox63}.  An elementary arrangement of plane
surfaces that elucidates the \textit{local} element of this
phenomenon is found upon considering the domain of $\mathbb{R}^3$
that is the union of the upper half-space together with all points
$(x,y,t)$ with $(x,y)$ in the first quadrant of the plane, i.e. the
union of a half-space and an infinite wedge.  The boundary consists
of $3$ faces:  the $4$th quadrants of both the $xt$ and $yt$-planes
and the piece of the $xy$-plane outside of the $1$st quadrant of the
$xy$-plane.  The requirement of postulate (iii) is met only by the
\textit{negative} $t$-axis.  But no color change is possible there
because any color on either of the $4$th quadrants must be continued
across the \textit{positive} $x$ or $y$-axis to the $3$rd face of
the boundary.  On the other hand, a color change is possible for the
complementary domain and is possible for the exterior domain to the
compound of $5$ tetrahedra.

Is there a polyhedral surface with a finite number of faces for
which both interior and exterior mixed problems are monochromatic?

\bibliographystyle{amsalpha}
\bibliography{neu2003}

\end{document}